\documentclass[reqno]{amsart}
\title[Congruences for two-color partitions]{Congruences for two-color partitions with odd smallest part}
\usepackage{amssymb,amsmath,amsthm,epsfig,graphics,latexsym,hyperref}
\theoremstyle{definition}
\newtheorem{definition}{Definition}
\theoremstyle{plain}
\newtheorem{lemma}      {Lemma}

\newtheorem{theorem}    {Theorem}

\newtheorem{corollary}  {Corollary}
\newtheorem{conjecture} {Conjecture}
\theoremstyle{remark}

\newtheorem{remark}{{\bf Remark}}
\numberwithin{equation}{section}
 \newcommand{\Mod}[1]{\ (\mathrm{mod}\ #1)}

\newcommand{\fr}{\frac}

\newmuskip\pFqskip
\mathchardef\pFcomma=\mathcode`, 
\newcommand*\pFq[5]{%
  \begingroup
  \begingroup\lccode`~=`,
  \lowercase{\endgroup\def~}{\pFcomma\mkern\pFqskip}%
 \mathcode`,=\string"8000
 {}_{#1}\phi_{#2}\biggl[\genfrac..{0pt}{}{#3}{#4};#5\biggr]%
 \endgroup
 }
 \newmuskip\pGqskip
\mathchardef\pGcomma=\mathcode`, 

\begin{document}
\author[ G. E. Andrews and M. El Bachraoui]{George E. Andrews and Mohamed El Bachraoui}
\address{The Pennsylvania State University, University Park, Pennsylvania 16802}
\email{andrews@math.psu.edu}
\address{Dept. Math. Sci,
United Arab Emirates University, PO Box 15551, Al-Ain, UAE}
\email{melbachraoui@uaeu.ac.ae}
\begin{abstract}
For a fixed positive integer $k$, let $C(k,n)$ denote the number of two-color partitions of $n$ with odd smallest part and restrictions on even parts, and let $C_k(q)$ be its generating function. We show that $C(1,n)\equiv d(2n-1)\pmod{4}$ and obtain congruences modulo $2$ and $4$ for $C(k,n)$ when $k=2,3$. Using $q$-series methods we derive closed formulas for $C_k(q)$ in terms of eta-quotients and formulate Ramanujan-type congruences for the limiting sequence arising from $\lim_{k\to\infty} C_k(q)$.
\end{abstract}
\date{\textit{\today}}
\thanks{First author partially supported by Simons Foundation Grant 633284}
\maketitle
\section{Introduction}
Throughout, we use the standard $q$-shifted factorial notation (see, for example, \cite{Andrews, Gasper-Rahman}):
\[
(a;q)_0 = 1,\  (a;q)_n = \prod_{j=0}^{n-1} (1-aq^j),\quad
(a;q)_{\infty} = \prod_{j=0}^{\infty} (1-aq^j),
\]
\[
(a_1,\ldots,a_k;q)_n = \prod_{j=1}^k (a_j;q)_n,\
(a_1,\ldots,a_k;q)_{\infty} = \prod_{j=1}^k (a_j;q)_{\infty}.
\]
Integer partitions in which each part may occur in two colors have been studied
extensively in recent years. For references on two-color partitions that are closely
related to the present work, see for
instance~\cite{Andrews 1987, Andrews 2021, Andrews-Bachraoui evens-1-color, Andrews-Bachraoui positive, Andrews-Kumar, Banerjee-et-al, Chern-Wang}.
In particular,
Banerjee, Bringmann, and Keith~\cite{Banerjee-et-al} have recently investigated a weighted
companion of the generating functions $C_k(q)$ considered here, providing
further evidence of the richness of these two-color partition families.
In this paper we consider sequences of integer partitions in two colors (blue and
red), as in the following definition.
\begin{definition}
Let $k$ be a fixed positive integer. For a positive integer $n$, let $\mathcal{C}(k,n)$ denote the set
of two-color partitions $\pi(n)$ of $n$ in which:
\begin{itemize}
\item the smallest part $s(\pi)$ is odd and occurs at least once in blue,
\item every even blue part is at least $2k-1$ greater than $s(\pi)$,
\item the even parts of the same color are distinct.
\end{itemize}
Accordingly, let $C(k,n)=|\mathcal{C}(k,n)|$.
Letting $q^{2n+1}\fr{(-q^{2n+2k};q^2)_\infty}{(q^{2n+1};q^2)_\infty}$ generate the blue parts and
$\fr{(-q^{2n+2};q^2)_\infty}{(q^{2n+1};q^2)_\infty}$ generate the red parts, we obtain
\begin{equation}\label{gen C}
C_{k}(q):= \sum_{n\geq 0}C(k,n) q^n
=\sum_{n\geq 0} q^{2n+1} \fr{(-q^{2n+2k},-q^{2n+2};q^2)_\infty}{(q^{2n+1};q^2)_\infty^2}.
\end{equation}
\end{definition}
For example,
\begin{align*}
C_1(q)&=q + 2 q^2 + 6 q^3 + 10 q^4 + 19 q^5 + 34 q^6 + 58 q^7 + 92 q^8 +150 q^9 +\cdots \\
C_2(q)&=q + 2 q^2 + 5 q^3 + 8 q^4 + 15 q^5 + 26 q^6 + 43 q^7 + 68 q^8 + 109 q^9 +\cdots \\
C_3(q)&=q + 2 q^2 + 5 q^3 + 8 q^4 + 14 q^5 + 24 q^6 + 39 q^7 + 60 q^8 + 95 q^9 +\cdots.
\end{align*}
We note that the previous definition is the case $m=1$ of~\cite[Definition~2.1]{Andrews-Bachraoui positive}.

Our first main goal in this paper is to establish the following congruences modulo $4$ for the coefficients of $C_k(q)$ for $k=1,2,3$.
The first congruence is expressed in terms of the function $d(N)$, counting the number of positive divisors of $N$.
\begin{theorem}\label{thm cong1}
For any positive integer $n$, we have
\[
C(1,n) \equiv d(2n-1) \Mod{4}.
\]
\end{theorem}
\begin{corollary}\label{cor odd-C1}
For any positive integer $n$, the number $C(1,n)$ is odd if and only if $2n-1$ is a perfect square. Equivalently,
the number $C(1,n)$ is odd if and only if $n=2k(k+1)+1$ for some nonnegative integer $k$.
\end{corollary}
\begin{corollary}\label{cor C1-mod4}
For $n\geq 1$, write $2n-1$ in the form
\[
2n-1 = \prod_{i=1}^r p_i^{e_i},
\]
where $p_1,\ldots,p_r$ are distinct odd primes and $e_i\geq 1$. Let $t$ be the number of indices $i$ with $e_i$ odd. Then
\[
C(1,n) \equiv
\begin{cases}
1\ \text{or}\ 3 \Mod{4}, & t=0, \\
2 \Mod{4}, & t=1, \\
0 \Mod{4}, & t\geq 2.
\end{cases}
\]
In particular, $C(1,n)\not\equiv 0\Mod{4}$ if and only if\  $2n-1$ is a square or a prime times a square.
\end{corollary}
\begin{corollary}\label{cor 9n+8}
For all nonnegative integers $n$ we have
\[
C(1,9n+8) \equiv 0 \Mod{4}.
\]
\end{corollary}
\medskip

Proofs of Corollaries~\ref{cor odd-C1}, \ref{cor C1-mod4}, and~\ref{cor 9n+8} are given in Section~\ref{sec proof 3-Corollaries}.
\begin{theorem}\label{thm cong2}
For any positive integer $n$, we have
\begin{align}
C(2,4n) &\equiv 0 \Mod{4}, \label{cong2-1} \\
C(2,4n-2) &\equiv 2 \Mod{4}. \label{cong2-2}
\end{align}
\end{theorem}
\begin{theorem}\label{thm cong3}
For any positive integer $n$, we have
\[
C(3,4n) \equiv 0 \pmod{4}.
\]
\end{theorem}

Our second goal in this paper is to prove the following formulas for $C_k(q)$ for $k=1,2,3$.
\begin{theorem}\label{thm C1}
There holds
\begin{equation}
C_1(q)=
\Big(
2q\sum_{n\ge0}
     \frac{q^{n^{2}+n}\bigl(1-q^{2n+1}\bigr)}
          {1+q^{2n+1}}
\Big) \fr{(-q^2;q^2)_\infty^2}{(q;q^2)_\infty^2}
-q\sum_{n\ge0}\frac{(-q)^{n}}{1+q^{2n+1}}.
\end{equation}
\end{theorem}
\begin{theorem}\label{thm C21}
There holds
\begin{equation}\label{C21-id}
C_{2}(q) = \fr{2q}{(1+q)^2}\fr{(-q^2;q^2)_\infty^2}{(q;q^2)_\infty^2} - \fr{q}{(1+q)^2}.
\end{equation}
\end{theorem}
\begin{theorem}\label{thm C31}
There holds
\begin{equation}\label{C31-id}
C_{3}(q) = \frac{2q(2q^2 - q + 1)}{(1+q)(1+q^2)(1+q^3)^2}\fr{(-q^2;q^2)_\infty^2}{(q;q^2)_\infty^2}
- \fr{q(1-q+q^2+q^3)}{(1+q)(1+q^3)^2}.
\end{equation}
\end{theorem}
Furthermore, we get the following immediate consequence of Theorem~\ref{thm C21}.
\begin{corollary}\label{cor cong2}
For any positive integer $n$, we have
\[
C(2,n) \equiv n \Mod{2}.
\]
\end{corollary}
Finally, noting that
\begin{align*}
\fr{q(1-q+q^2+q^3)}{(1+q)(1+q^3)^2} &\equiv \fr{q(1+q+q^2+q^3)}{(1+q)(1+q^3)^2} \Mod{2} \\
&= \fr{q(1+q^2)}{(1+q^3)^2} \Mod{2} \\
&= \sum_{n\geq 0} (q^{6n+1}+q^{6n+3}) \Mod{2},
\end{align*}
and employing Theorem~\ref{thm C31}, we get the following congruence.
\begin{corollary}\label{cor cong3}
For any positive integer $n$, we have
\begin{equation*}
C(3,n) \equiv 1\Mod{2}\ \text{if and only if\ } n=1, 3 \Mod{6}.
\end{equation*}
\end{corollary}
%
Beyond their combinatorial definition, the generating functions $C_k(q)$ also admit a modular-forms interpretation; see, for example, Ono~\cite{Ono} for background on modular forms and $q$-series. In particular, the common factor
\[
\frac{(-q^2;q^2)_\infty^2}{(q;q^2)_\infty^2}
\]
appearing in our formulas for $C_1(q)$, $C_2(q)$, and $C_3(q)$ can be written as an eta-quotient of weight~$0$ on $\Gamma_0(4)$. We return to this point in Section~\ref{sec remarks} (see Remark~\ref{rem eta-quotient}).
\section{$q$-series background}\label{sec background}
In this section we collect several well-known facts from the theory of basic hypergeometric series that will be needed in our proofs.
Recall the definition of the basic hypergeometric series~\cite[Eq.~(1.2.22)]{Gasper-Rahman}:
\[
  \pFq{\,r+1\,}{\,r\,}
      {a_1, a_2, \dots, a_{r+1}}
      {b_1, b_2, \dots, b_r}
      {q; z}
  =
  \sum_{n=0}^{\infty}
    \frac{(a_1;q)_n (a_2;q)_n \cdots (a_{r+1};q)_n}
         {(q;q)_n (b_1;q)_n \cdots (b_r;q)_n}
    \, z^n.
\]

Heine's first transformation~\cite[Corollary~2.3]{Andrews} states that
\begin{equation}\label{Heine-1}
\pFq{2}{1}{a,b}{c}{q, z}
= \fr{(b,az;q)_\infty}{(c,z;q)_\infty}
  \pFq{2}{1}{c/b,z}{az}{q, b}.
\end{equation}
The Rogers–Fine identity~\cite[Eq.~(2.7)]{Andrews 1981} states that
\begin{equation}\label{Rogers-Fine}
\sum_{n\geq 0}\fr{(\alpha;q)_n \tau^n}{(\beta;q)_n}
=\sum_{n\geq 0}
  \fr{(\alpha;q)_n (q\alpha\tau/\beta;q)_n \beta^n \tau^n q^{n^2-n}
       (1-\alpha\tau q^{2n})}
     {(\beta;q)_n (\tau;q)_{n+1}}.
\end{equation}

We shall also make use of Watson's transformation~\cite[Eq.~(III.17)]{Gasper-Rahman},
which in our notation reads
\[
  \pFq{8}{7}
      {a, qa^{1/2}, -qa^{1/2}, b, c, d, e, f}
      {a^{1/2}, -a^{1/2}, aq/b, aq/c, aq/d, aq/e, aq/f}
      {q, a^{2}q^{2}/(bcdef)}
  =
\]
\begin{equation}\label{Watson}
  \frac{(aq, aq/de, aq/df, aq/ef; q)_\infty}
       {(aq/d, aq/e, aq/f, aq/def; q)_\infty}
  \,
  \pFq{4}{3}
      {aq/bc, d, e, f}
      {aq/b, aq/c, def/a}
      {q, q},
\end{equation}
as well as the following three-term transformation formula~\cite[Eq.~(III.34)]{Gasper-Rahman}:
\[
  \pFq{3}{2}
      {a,b,c}
      {d,e}
      {q;\dfrac{de}{abc}}
  =
  \frac{(e/b;q)_\infty (e/c;q)_\infty}{(e;q)_\infty (e/bc;q)_\infty}
  \,\pFq{3}{2}
      {d/a,b,c}
      {d,bcq/e}
      {q;q}
\]
\begin{equation}\label{three-term}
 +
  \frac{(d/a;q)_\infty (b;q)_\infty (c;q)_\infty (de/bc;q)_\infty}
       {(d;q)_\infty (e;q)_\infty (bc/e;q)_\infty (de/abc;q)_\infty}
  \,\pFq{3}{2}
      {e/b,e/c,de/abc}
      {de/bc,eq/bc}
      {q;q}.
\end{equation}

We also need the following identity~\cite{Andrews-Bachraoui, Chan-Mao}:
\begin{equation}\label{Chan-Mao}
\sum_{n\geq 0}\fr{(y,q/y;q)_n q^n}{(z,q/z;q)_{n+1}}
= C_q(y,z)
  +\fr{(y,q/y;q)_\infty}{z(1-y/z)(1-q/yz)(z,q/z;q)_\infty},
\end{equation}
where
\[
C_q(y,z) = \fr{1}{y\big(1-\frac{z}{y}\big) \big(1-\frac{q}{yz}\big)}.
\]
Note that
\begin{equation}\label{basic-C}
C_{q^2}(q,-q^{-2})= \fr{q^2}{(1+q^3)^2}
\quad\text{and}\quad
C_{q^2}(q,-1)=\fr{1}{(1+q)^2}.
\end{equation}

We will also require the following formula from~\cite{Andrews-Bachraoui companions}:
\begin{equation}\label{Gea-Mel cmpn}
\sum_{n\geq 1} \fr{q^n (qy,q/y;q)_{n-1}}{(qz,q/z;q)_{n}}
= \fr{1}{y+y^{-1}-z-z^{-1}}
  \Big( 1 - \fr{(qy,q/y;q)_\infty}{(qz,q/z;q)_\infty}\Big)
\end{equation}
and the following identity from~\cite{Andrews-Bachraoui}:
\begin{align}
\sum_{n\geq 0} \fr{(y,q/y;q)_n q^{2n+1}}{(qz,q/z;q)_{n+1}}
&=(1-q)C_q(y,z)C_q(y,qz)
   \fr{(y,q/y;q)_\infty}{(qz,q/z;q)_\infty} \nonumber \\[4pt]
&\quad + \fr{(1-q)(1-z) C_q(y,z)C_q(y,qz)}{z}
       + C_q(y,z). \label{id appl-1}
\end{align}

Finally, we will use the following identity of Andrews~\cite[Eq.~(5.1)]{Andrews 1979}:
\begin{equation}\label{marvelous}
\sum_{n\geq 0} \fr{q^{n}}{(-q;q^2)_{n+1}}
 = \sum_{n\geq 0} (-1)^n q^{2n(3n+2)}(1+q^{4n+2})
\end{equation}
and the following identity of Ramanujan~\cite[Entry~9.5.1]{Andrews-Berndt I}:
\begin{equation}\label{Entry951}
\sum_{n\geq 0}\fr{(-1)^n a^{2n} q^{n(n+1)}}{(-aq;q^2)_{n+1}}
= \sum_{n\geq 0} a^{3n} q^{3n^2 +2n}(1-aq^{2n+1}).
\end{equation}
\section{Proof of Theorem~\ref{thm cong1}}\label{sec proof cong1}
We need two lemmas.
\begin{lemma}\label{lem1 cong1}
There holds
\[
\sum_{n\geq 1}\fr{q^n}{1-q^{2n-1}} = \sum_{n\geq 1} d(2n-1)q^n.
\]
\end{lemma}
\begin{proof}
We first note the basic fact
\[
\sum_{n\geq 1}d(n) q^{n} = \sum_{n\geq 1}\fr{q^n}{1-q^n}.
\]
Hence,
\begin{align*}
\sum_{n\geq 1} d(2n-1)q^{2n-1} &= \fr{1}{2} \sum_{n\geq 1}\Big(\fr{q^n}{1-q^n}-\fr{(-q)^n}{1-(-q)^n}\Big) \\
&=\fr{1}{2} \sum_{n\geq 1}\Big(\fr{q^{2n-1}}{1-q^{2n-1}}+\fr{q^{2n-1}}{1+q^{2n-1}}\Big) \\
&=\fr{1}{2} \sum_{n\geq 1} \fr{2q^{2n-1}}{1-q^{4n-2}} \\
&=\sum_{n\geq 1}\fr{q^{2n-1}}{1-q^{4n-2}}.
\end{align*}
Now multiply both sides by $q$ and then replace $q^2$ by $q$ to get the desired result.
\end{proof}
For a power series $A(q):= \sum_{n\geq 0}a_n q^n$
and a nonnegative integer $N$, we let $[q^N]A(q):= a_n$.
Our next lemma is in terms of the set $\mathbb{T}$ given by
\[
\mathbb{T} =\{ n\in\mathbb{N}_0\ :\ n=k(k+1)\ \text{for some $k\in\mathbb{N}_0$} \}.
\]
\begin{lemma}\label{lem2 cong1}
For any nonnegative integer $n$, we have
\begin{align*}
[q^n]\fr{(-q^2;q^2)_\infty^2}{(q;q^2)_\infty^2} &\equiv [q^n]\fr{(q^2;q^2)_\infty^2}{(q;q^2)_\infty^2} \Mod{2} \\
&\equiv 1\Mod{2}\ \text{if and only if $n\in\mathbb{T}$.}
\end{align*}
\end{lemma}
\begin{proof}
The first congruence in~\eqref{lem2 cong1} follows from the following basic congruence
\begin{equation}\label{basic-cong}
(1-q)^2\equiv (1+q)^2\Mod{4}.
\end{equation}
For the second congruence, recall first the Gauss formula~\cite[(2.2.13)]{Andrews}
\begin{equation}\label{Gauss}
\fr{(q^2;q^2)_\infty}{(q;q^2)_\infty} = \sum_{n\geq 0} q^{\fr{n(n+1)}{2}}.
\end{equation}
Then we get
\begin{align*}
\fr{(q^2;q^2)_\infty^2}{(q;q^2)_\infty^2} &= \Big(\sum_{n\geq 0} q^{\fr{n(n+1)}{2}} \Big)^2 \\
&\equiv \sum_{n\geq 0} \big(q^{\fr{n(n+1)}{2}} \big)^2 \Mod{2} \\
&\equiv \sum_{n\geq 0} q^{n(n+1)} \Mod{2},
\end{align*}
implying the desired congruence.
\end{proof}
\emph{First proof of Theorem~\ref{thm cong1}.\ }
Note that
\[
\fr{q^{n^2+n}(1-q^{2n+1})}{1+q^{2n+1}} = q^{n^2+n} \Big( 1-\fr{2 q^{2n+1}}{1+q^{2n+1}} \Big),
\]
which combined with Theorem~\ref{thm C1} yields
\begin{align}
C_1(q) &= \Big(2q\sum_{n\geq 0} \fr{q^{n^{2}+n}\bigl(1-q^{2n+1}\bigr)}{1+q^{2n+1}}\Big) \fr{(-q^2;q^2)_\infty^2}{(q;q^2)_\infty^2}
-q\sum_{n\geq 0}\fr{(-q)^{n}}{1+q^{2n+1}} \nonumber \\
&=\left(2q\sum_{n\geq 0} q^{n^{2}+n} \Big(1-2\fr{q^{2n+1}}{1+q^{2n+1}}\Big) \right)\fr{(-q^2;q^2)_\infty^2}{(q;q^2)_\infty^2}
-q\sum_{n\geq 0}\fr{(-q)^{n}}{1+q^{2n+1}} \nonumber \\
&\equiv 2q\fr{(q^2;q^2)_\infty^2}{(q;q^2)_\infty^2} \sum_{n\geq 0} q^{n^{2}+n} - q\sum_{n\geq 0}\fr{(-q)^{n}}{1+q^{2n+1}} \Mod{4}. \label{help1 cong1}
\end{align}
Note that by Lemma~\ref{lem2 cong1}, if $n\not\in\mathbb{T}$, then
\begin{equation}\label{help2 cong1}
[q^n]\fr{(q^2;q^2)_\infty^2}{(q;q^2)_\infty^2} \equiv 0\Mod{2}.
\end{equation}
On the other hand, if $n\in\mathbb{T}$, say $n=k(k+1)$ for some $k\in\mathbb{N}_0$, then
\begin{equation}\label{help3 cong1}
[q^n] \Big(\sum_{m\geq 0} q^{m^2+m+1} \Big) \equiv 0\Mod{2}
\end{equation}
because the equation $k^2+k = m^2+m+1$ has no solutions in $\mathbb{N}_0$.
Then appealing to the identities~\eqref{help1 cong1}-\eqref{help3 cong1}, we achieve
\begin{align*}
C_1(q) &\equiv -q\sum_{n\geq 0}\fr{(-q)^{n}}{1+q^{2n+1}} \Mod{4} \\
& = -\sum_{n\geq 1} \fr{(-q)^n}{1+q^{2n-1}} \Mod{4}.
\end{align*}
Equivalently,
\begin{equation}\label{help4 cong1}
C_1(-q) \equiv -\sum_{n\geq 1} \fr{q^n}{1-q^{2n-1}} \Mod{4}.
\end{equation}
Besides, from~\eqref{gen C} and~\eqref{basic-cong},
we derive
\begin{align}
C_1(-q) &= -\sum_{n\geq 0}\fr{q^{2n+1} (-q^{2n+2};q^2)_\infty^2}{(-q^{2n+1};q^2)_\infty^2}  \nonumber \\
&\equiv -\sum_{n\geq 0}\fr{q^{2n+1} (-q^{2n+2};q^2)_\infty^2}{(q^{2n+1};q^2)_\infty^2} \Mod{4} \nonumber \\
&= -C_1(q) \Mod{4}. \label{help5 cong1}
\end{align}
Now combine~\eqref{help5 cong1} and~\eqref{help4 cong1} with Lemma~\ref{lem1 cong1} to obtain
\begin{align*}
C_1(q)&\equiv \sum_{n\geq 1} \fr{q^n}{1-q^{2n-1}} \Mod{4} \\
&= \sum_{n\geq 1} d(2n-1)q^n,
\end{align*}
which is the desired result.

\emph{Second proof of Theorem~\ref{thm cong1}.\ } By virtue of~\eqref{basic-cong},~\eqref{Heine-1}, and Lemma~\ref{lem1 cong1}, we get
\begin{align*}
\sum_{n\geq 0} C(1,n) q^n &=\sum_{n\geq 0} \fr{q^{2n+1}(-q^{2n+2};q^2)_\infty^2}{(q^{2n+1};q^2)_\infty^2} \\
&= \sum_{n\geq 0} \fr{q^{2n+1}(q^{2n+2};q^2)_\infty^2}{(q^{2n+1};q^2)_\infty^2} \Mod{4} \\
&=q\fr{(q^2;q^2)_\infty^2}{(q;q^2)_\infty^2} \pFq{2}{1}{q,q}{q^{2}}{q^2, q^{2}} \\
&=q\fr{(q^2;q^2)_\infty^2}{(q;q^2)_\infty^2} \fr{(q,q^3;q^2)_\infty}{(q^2;q^2)_\infty^2} \pFq{2}{1}{q,q^2}{q^{3}}{q^2, q} \\
&= \fr{q}{1-q}\sum_{n\geq 0} \fr{q^n (q;q^2)_n}{(q^3;q^2)_n} \\
&=\sum_{n\geq 0} \fr{q^{n+1} (q;q^2)_n}{(q;q^2)_{n+1}} \\
&=\sum_{n\geq 1} \fr{q^{n}}{1-q^{2n-1}} \\
&=\sum_{n\geq 1} d(2n-1)q^n,
\end{align*}
as desired.
\section{Proofs of Corollaries~\ref{cor odd-C1}-\ref{cor 9n+8}}\label{sec proof 3-Corollaries}
%
We start with a lemma.
\begin{lemma}\label{lem d-mod4}
Let $N$ be a positive odd integer with prime factorisation
\[
N = \prod_{i=1}^{r} p_i^{e_i},
\]
where $p_1,\ldots,p_r$ are distinct odd primes and $e_i\geq 1$. Let $t$ be the number of indices $i$ for which $e_i$ is odd. Then
\[
d(N) \equiv
\begin{cases}
1\ \text{or}\ 3 \Mod{4}, & t=0, \\
2 \Mod{4}, & t=1, \\
0 \Mod{4}, & t\geq 2.
\end{cases}
\]
\end{lemma}
\begin{proof}
We have
\begin{align*}
d(N) &= \prod_{i=1}^{r} (e_i+1) \\
&= \prod_{e_i\ \text{odd}}(e_i+1) \prod_{e_j\ \text{even}} (e_i+1) \\
&=\prod_{e_i\ \text{odd}} (2a_i+2) \prod_{e_j\ \text{even}} (2b_j+1) \\
&=2^{t} M,
\end{align*}
where $M$ is an odd integer. If $t=0$, then $d(N)=M$ is odd, hence $d(N)\equiv 1$ or $3\Mod{4}$. If $t=1$, then $d(N)=2M\equiv 2\Mod{4}$, and if $t\geq 2$, then $4\mid d(N)$, so $d(N)\equiv 0\Mod{4}$. This confirms the desired congruences.
\end{proof}

\emph{Proof of Corollary~\ref{cor odd-C1}.\ }
By Theorem~\ref{thm cong1} we have
\[
C(1,n)\equiv d(2n-1)\Mod{4}
\]
for all positive integers $n$. In particular, $C(1,n)$ is odd if and only if $d(2n-1)$ is odd. Writing
\[
2n-1 = \prod_{i=1}^{r} p_i^{e_i}
\]
as in Lemma~\ref{lem d-mod4}, we see that $d(2n-1)$ is odd if and only if $t=0$, i.e., all exponents $e_i$ are even. This is equivalent to $2n-1$ being a perfect square.

Since $2n-1$ is odd, its square root $s$ must be odd, say $s=2k+1$ with $k\in\mathbb{N}_0$. Then
\[
2n-1 = (2k+1)^2 = 4k^2+4k+1,
\]
so $n=2k(k+1)+1$, as desired.

\emph{Proof of Corollary~\ref{cor C1-mod4}.\ }
Again by Theorem~\ref{thm cong1} we have $C(1,n)\equiv d(2n-1)\Mod{4}$. Writing $2n-1$ as in the statement of the corollary and applying Lemma~\ref{lem d-mod4} with $N=2n-1$ gives the three cases for $C(1,n)\Mod{4}$.

For the final assertion, observe that $C(1,n)\not\equiv 0\Mod{4}$ if and only if $t\leq 1$. If $t=0$, then all exponents $e_i$ are even and $2n-1$ is a square. If $t=1$, say $e_1$ is odd and $e_i$ is even for $i\geq 2$, then we may write
\[
2n-1 = p_1^{2a+1} \prod_{i=2}^{r} p_i^{2b_i} = p_1\Big( p_1^{a} \prod_{i=2}^{r} p_i^{b_i} \Big)^2,
\]
so $2n-1$ is a prime times a square. This completes the proof.

\emph{Proof of Corollary~\ref{cor 9n+8}.\ }
By Theorem~\ref{thm cong1} we have
\[
C(1,9n+8)\equiv d(18n+15)\Mod{4}
\]
for all $n\geq 0$. Since $18n+15 = 3(6n+5)$, the prime $3$ occurs to an odd exponent in the factorisation of $18n+15$.

We claim that $6n+5$ is never a perfect square. Indeed, if $6n+5 = s^2$ for some integer $s$, then reducing modulo $6$ would give $s^2\equiv 5\Mod{6}$, which is impossible because a square modulo $6$ can only be $0,1,3,$ or $4$. Thus $6n+5$ is not a square, so in the factorisation of $18n+15$ at least one other prime besides $3$ must also occur to an odd exponent.

In the notation of Lemma~\ref{lem d-mod4} we therefore have $t\geq 2$, and hence $d(18n+15)\equiv 0\Mod{4}$. This gives $C(1,9n+8)\equiv 0\Mod{4}$, as required.
\section{Proof of Theorem~\ref{thm cong2}}\label{sec proof cong2}
We start by a lemma.
\begin{lemma}\label{lem1 cong2}
Let
\[
\overline{C}_2(q) := \fr{q}{1+q^2}\fr{(q^2;q^2)_\infty^2}{(q;q^2)_\infty^2}
\pFq{3}{2}{-q^2,q,q}{q^{2},-q^{4}}{q^{2};q^2}.
\]
Then
\begin{equation}\label{C2-}
\overline{C}_2(q)
=\fr{q(+q)}{(1+q^2)(1-q)}+2\sum_{n\geq 1} \fr{q^{n^2+3n+1}}{(1+q^{2n})(1+q^{2n+2})}\Big(1+\fr{2q^{2n+1}}{1-q^{2n+1}}\Big).
\end{equation}
\end{lemma}
\begin{proof}
An application of~\eqref{Watson} with $a=q$ and $f\to\infty$ gives after simplification
\[
\sum_{n\geq 0} \fr{(1-q^{2n+1})}{(1-q)}\fr{(b,c,d,e;q)_n}{(q^2/b,q^2/c,q^2/d,q^2/e;q)_n} (-1)^n q^{{n\choose 2}}\fr{q^{4n}}{(bcde)^n}
\]
\[
=\fr{(q^2,q^2/de;q)_\infty}{(q^2/d,q^2/e;q)_\infty}\pFq{3}{2}{q^2/bc,d,e}{q^{2}/b,q^{2}/c}{q;q^2/de}.
\]
Now, in the foregoing identity, let $q\to q^2$ and $(b,c,d,e)=(q^2,-1,q,q)$ and rearrange to obtain
\[
\fr{(q^3,q^3;q^2)_\infty}{(q^4,q^2;q^2)_\infty}\sum_{n\geq 0} \fr{(1-q^{4n+2})}{(1-q^2)}\fr{(q^2,-1,q,q;q^2)_n}{(q^2,-q^4,q^3,q^3;q^2)_n} (-1)^n q^{n^2-n}\fr{q^{8n}}{(-q^4)^n}
\]
\[
=\pFq{3}{2}{-q^2,q,q}{q^{2},-q^{4}}{q^{2};q^2}.
\]
Multiplying both sides by $\fr{q}{1+q^2}\fr{(q^2;q^2)_\infty^2}{(q;q^2)_\infty^2}$ and simplifying yield
\begin{align*}
\overline{C_2}(q) &=\fr{q(1-q^2)}{(1+q^2)(1-q)^2} \sum_{n\geq 0}\fr{(1-q^{4n+2})}{(1-q^2)} \fr{(-1,q,q;q^2)_n}{(-q^4,q^3,q^3;q^2)_n} q^{n^2+3n} \\
&=\fr{q}{1+q^2} \sum_{n\geq 0} \fr{(1+q^{2n+1})}{(1-q^{2n+1})}\fr{2(1+q^2)}{(1+q^{2n})(1+q^{2n+2})} q^{n^2+3n} \\
&=\fr{q}{1+q}{(1+q^2)(1-q)}+ 2\sum_{n\geq 1} \fr{q^{n^2+3n+1}}{(1+q^{2n})(1+q^{2n+2})}\Big(1+\fr{2q^{2n+1}}{1-q^{2n+1}}\Big),
\end{align*}
as desired.
\end{proof}
By~\eqref{basic-cong} and Lemma~\ref{lem1 cong2}, we have
\begin{align}
C_2(q)&= \sum_{n\geq 0}\fr{q^{2n+1}(-q^{2n+2},-q^{2n+4};q^2)_\infty}{(q^{2n+1};q^2)_\infty^2} \nonumber \\
&\equiv \sum_{n\geq 0}\fr{q^{2n+1}(q^{2n+2},q^{2n+4};q^2)_\infty}{(q^{2n+1};q^2)_\infty^2} \Mod{4} \nonumber \\
&=\fr{q}{1+q^2}\fr{(q^2;q^2)_\infty^2}{(q;q^2)_\infty^2} \pFq{3}{2}{-q^2,q,q}{q^{2},-q^{4}}{q^{2};q^2} \nonumber \\
&=\overline{C}_2(q) \\
&=\fr{q(+q)}{(1+q^2)(1-q)}+2\sum_{n\geq 1} \fr{q^{n^2+3n+1}}{(1+q^{2n})(1+q^{2n+2})}\Big(1+\fr{2q^{2n+1}}{1-q^{2n+1}}\Big) \nonumber\\
&\equiv \fr{q(+q)}{(1+q^2)(1-q)} + 2\sum_{n\geq 1}\fr{q^{n^2+3n+1}}{(1+q^{2n})(1+q^{2n+2})} \Mod{4}. \label{help cong2}
\end{align}
Note that the second term of the right-hand side of~\eqref{help cong2} is an odd function of $q$ and its first term is
\[
\fr{q(1+q)}{(1+q^2)(1-q)} = \fr{q(1+q)^2}{1-q^4} = \fr{q+2q^2+q^3}{1-q^4}.
\]
Consequently, combining the above notes with the fact that $C_2(q)\equiv \overline{C}_2(q) \Mod{4}$, see that
\[
C(2,4n)\equiv 0\Mod{4}\ \text{and\ } C(2,4n+2)\equiv 2\Mod{4}.
\]
This completes the proof.
\section{Proof of Theorem~\ref{thm cong3}}\label{sec proof cong3}
Throughout we will be using the following basic congruences
\begin{equation}\label{basic-cong mod2}
1+q \equiv 1-q \Mod{2}\ \text{and\ }
(1+q)^2 \equiv 1+q^2\Mod{2},
\end{equation}
We need a lemma.
\begin{lemma}\label{lem cong3}
Let
\[
f(q)=\frac{q^5(1-q^2)^2(1-q^4)}{(1+q^2)(1+q^4)(1-q)^2(1-q^3)^2}\ \text{and\ }
g(q)=\frac{q^5(1-q^4)}{(1+q^2)(1-q)^2(1-q^3)^2}.
\]
Then for all integers $n\ge0$ we have
\[
[q^{4n}]f(q)\equiv [q^{4n}]g(q)\Mod{4}.
\]
\end{lemma}
\begin{proof}
We compute
\[
\begin{aligned}
f(q)-g(q)
&=\frac{q^5(1-q^4)\bigl((1-q^2)^2-(1+q^4)\bigr)}
        {(1+q^2)(1+q^4)(1-q)^2(1-q^3)^2} \\
&=\frac{-2q^7(1-q^4)}{(1+q^2)(1+q^4)(1-q)^2(1-q^3)^2} \\
&=-2q^7 h(q),
\end{aligned}
\]
where
\[
h(q)=\fr{(1-q^4)}{(1+q^2)(1+q^4)(1-q)^2(1-q^3)^2} := \sum_{n\ge0} c_n q^n.
\]
Then
\[
[q^{4n}]\big( f(q)-g(q) \big)=-2\,c_{4n-7}.
\]
Thus to prove $[q^{4n}]f(q)\equiv [q^{4n}]g(q)\Mod{4}$, it suffices to show that
$c_{4n-7}$ is even for all $n$.
With the help of~\eqref{basic-cong mod2}, we get
\begin{align*}
h(q) &\equiv \frac{1+q^4}{(1+q^2)(1+q^4)(1+q)^2(1+q^3)^2} \Mod{2} \\
&\equiv \frac{1}{(1+q^2)(1+q)^2(1+q^3)^2}\Mod{2} \\
&\equiv \frac{1}{(1+q^2)(1+q^2)(1+q^6)} \Mod{2}\\
&\equiv \frac{1}{(1+q^2)^2(1+q^6)}\Mod{2} \\
&\equiv \Big(\sum_{k\ge0} q^{2k}\Big)^2 \Big(\sum_{l\ge0} q^{6l}\Big) \Mod{2} \\
&\equiv \Big(\sum_{k\ge0} q^{4k}\Big)\Big(\sum_{l\ge0} q^{6l}\Big) \\
&=\sum_{k,l\geq 0} q^{4k+6l}.
\end{align*}
Consequently, for any integer $n\geq 0$,
\[
[q^{2n+1}]h(q)\equiv 0 \Mod{2}
\]
and therefore
\[
c_{4n-7} \equiv [q^{4n-7}]h(q)\equiv 0\Mod{2}.
\]
This completes the proof.
\end{proof}
We are now ready to prove the theorem.
We have
\begin{align}
C_3(q)&= \sum_{n\geq 0}\fr{q^{2n+1}(-q^{2n+2},-q^{2n+6};q^2)_\infty}{(q^{2n+1};q^2)_\infty^2} \nonumber \\
&=\sum_{n\geq 0}\fr{q^{2n+1}(-q^{2n+2},-q^{2n+4};q^2)_\infty}{(1+q^{2n+4})(q^{2n+1};q^2)_\infty^2} \nonumber \\
&=\sum_{n\geq 0}\fr{q^{2n+1}(-q^{2n+2},-q^{2n+4};q^2)_\infty \big( (1+q^{2n+4})-q^{2n+4}}{(1+q^{2n+4})(q^{2n+1};q^2)_\infty^2} \nonumber \\
&= C_2(q) -\sum_{n\geq 0}\fr{q^{4n+5}(-q^{2n+2},-q^{2n+6};q^2)_\infty}{(q^{2n+1};q^2)_\infty^2}.
\end{align}
For simplicity, let
\[
S(q):=\sum_{n\geq 0}\fr{q^{4n+5}(-q^{2n+2},-q^{2n+6};q^2)_\infty}{(q^{2n+1};q^2)_\infty^2}.
\]
Now we know by Theorem~\ref{thm cong2} that the coefficients of $q^{4n}$ in $C_2(q)$ are divisible by $4$. Hence to prove
that they are also in $C_3(q)$, we need only prove that the coefficients of $q^{4n}$ in $S(q)$ are divisible by $4$.
Now,
\begin{equation}\label{S id1}
S(q) = \fr{q^5 (-q^2,-q^6;q^2)_\infty}{(q;q^2)_\infty^2}\pFq{3}{2}{q,q,q^2}{-q^{2},-q^{6}}{q^{2};q^4}.
\end{equation}
We may now apply Watson's transformation~\eqref{Watson} with $q\to q^2$, $f\to\infty$, and $(a,b,c,d,e)=(q^4,-q^4,-1,q,q)$.
Then the left hand-side of~\eqref{Watson} becomes
\begin{align}
&\sum_{n\geq 0}\fr{(q^4,q^4,-q^4,-q^4,-1,q,q;q^2)_n}{(q^2,q^2,-q^2,-q^2,-q^6,q^5,q^5;q^2)_n} q^{6n} \lim_{f\to\infty}\fr{(f;q^2)_n}{f^n} \nonumber \\
&=1 + \sum_{n\geq 1}\fr{(1-q^{2n+2})}{(1-q^2)}\fr{(1-q^{4n+4})}{(1-q^4)}\fr{(1+q^{2n+2})}{(1+q^2)}\fr{2(1+q^2)(1+q^4)}{(1+q^{2n})(1+q^{2n+2})(1+q^{2n+4})} \nonumber \\
&\qquad \qquad \cdot \fr{(1-q)^2 (1-q^3)^2}{(1-q^{2n+1})^2 (1-q^{2n+3})^2} (-1)^n q^{n^2+5n}.\label{LHS help}
\end{align}
Noting that by l'Hopital's rule,
\[
\lim_{f\to\infty}\fr{(f;q^2)_n}{(\fr{f}{q^2};q^2)_n} = q^{2n}
\]
we see that the right-hand side of~\eqref{Watson} becomes
\begin{equation}\label{RHS help}
\fr{(q^6,q^4;q^2)_\infty}{(q^5;q^2)_\infty^2} \lim_{f\to\infty}\pFq{4}{3}{q,q,q^2,f}{-q^{2},-q^{6},f/q^2}{q^{2};q^2}
=\fr{(q^6,q^4;q^2)_\infty}{(q^5;q^2)_\infty^2} \pFq{3}{2}{q,q,q^2}{-q^{2},-q^{6}}{q^{2};q^4}.
\end{equation}
Then combining~\eqref{S id1}-~\eqref{RHS help}, we obtain
\begin{align*}
S(q) &= \fr{(-q^2,-q^6;q^2)_\infty}{(q;q^2)_\infty^2}\fr{(q^5;q^2)_\infty^2}{(q^6,q^4;q^2)_\infty}
\Big( q^5+2 q^5 \fr{(1+q^2)(1+q^4)(1-q)^2(1-q^3)^2}{(1-q^2)(1-q^4)(1+q^2)} \\
&\qquad \cdot \sum_{n\geq 1}\fr{(-1)^n q^{n^2+5n}(1-q^{2n+2})(1-q^{4n+4})(1+q^{2n+2})}{(1+q^{2n})(1+q^{2n+2})(1+q^{2n+4})(1-q^{2n+1})^2 (1-q^{2n+3})^2} \Big) \\
&= R_1(q) + R_2(q).
\end{align*}
Using~\eqref{basic-cong}, we deduce that
\begin{align*}
R_1(q) &= \fr{q^5 (-q^2,-q^6;q^2)_\infty (q^5;q^2)_\infty^2}{(q;q^2)_\infty^2(q^6,q^4;q^2)_\infty} \\
&=\fr{q^5 (-q^2;q^2)_\infty^2}{(q^2;q^2)_\infty^2} \fr{(1-q^2)^2 (1-q^4)}{(1+q^2)(1+q^4)(1-q)^2 (1-q^3)^2} \\
&\equiv \fr{q^5 (1-q^2)^2 (1-q^4)}{(1+q^2)(1+q^4)(1-q)^2 (1-q^3)^2}.
\end{align*}
Next using~\eqref{basic-cong} and appealing to Lemma~\ref{lem cong3},
we get
\begin{align*}
[q^{4n}]R_1(q)
& \equiv [q^{4n}] \fr{q^5 (1-q^4)}{(1+q^2)(1-q)^2 (1-q^3)^2} \Mod{4}. \\
&=[q^{4n}] \fr{q^5(1-q^2)(1+q+q^2+q^3)^2 (1+q^3+q^6+q^9)^2}{(1-q^4)^2 (1-q^{12})^2} \\
&=[q^{4n}] \fr{1}{(1-q^4)^2 (1-q^{12})^2}
\Big(
     q^5 + 2q^6 + 2q^7 + 4q^8 + 4q^9 \\
& + 2q^{10} + 5q^{11}
     + 4q^{12} + q^{13} + 6q^{14} + 4q^{15} \\
 &+ 5q^{17} - 5q^{19} - 4q^{21}
    - 6q^{22} - q^{23} - 4q^{24} - 5q^{25} \\
&   - 2q^{26} - 4q^{27} - 4q^{28} - 2q^{29} - 2q^{30} - q^{31}
\Big).
\end{align*}
Note that all the coefficients in the foregoing expansion must be divisible by $4$ since this is true
in the last numerator and the denominator is a function of $q^4$.

Finally, using~\eqref{basic-cong mod2} we obtain
\begin{align*}
R_2(q) &= 2 q^5 \fr{(-q^2,-q^6;q^2)_\infty (q^5;q^2)_\infty^2 (1-q)^2(1-q^3)^3}{(q^6,q^4;q^2)_\infty (q;q^2)_\infty^2} \\
&\cdot \sum_{n\geq 1} \Big( \fr{(1-q^{2n+2})(1-q^{4n+4})(1+q^{2n+2})}{(1-q^2)(1-q^4)(1+q^2)} \fr{(1+q^2)(1+q^4)}{(1+q^{2n})(1+q^{2n+2})(1+q^{2n+4})} \\
& qquad \cdot\fr{(-1)^n q^{n^2+5n}}{(1-q^{2n+1})^2(1-q^{2n+2})^2} \Big)\\
&= 2q^5 \fr{(-q^2,-q^6;q^2)_\infty}{(q^6,q^4;q^2)_\infty} \\
&\sum_{n\geq 1}\fr{(1-q^{2n+2})(1-q^{4n+4})}{(1-q^2)(1-q^4)(1+q^2)}
\fr{(1+q^2)(1+q^4)}{(1+q^{2n})(1+q^{2n+4})}\fr{(-1)^n q^{n^2+5n}}{(1-q^{2n+1})^2(1-q^{2n+3})^2} \\
&\equiv 2q^5 \fr{(-q^2,-q^6;q^2)_\infty}{(q^6,q^4;q^2)_\infty} \\
&\sum_{n\geq 1}\fr{(1-q^{2n+2})(1-q^{4n+4})}{(1-q^2)(1-q^4)(1+q^2)}
\fr{(1+q^2)(1+q^4)}{(1+q^{2n})(1+q^{2n+4})}\fr{(-1)^n q^{n^2+5n}}{(1+q^{4n+2})(1+q^{4n+6})} \Mod{4} .
\end{align*}
But the last expression is an odd function of $q$ owing to the $q^5$ out front. Hence $R_2(q)$ has zero contribution to the
coefficients of the even powers of $q$. Therefore
\[
[q^{4n}] S(q) \equiv [q^{4n}] C_2(q) \equiv 0\Mod{4}.
\]
So, since
\[
C_3(q) = C_2(q)-S(q),
\]
we deduce that
\[
[q^{4n}] C_3(q) \equiv 0\Mod{4},
\]
which is the desired result.
\section{Proof of Theorem~\ref{thm C1}}\label{sec proof C1}
In the transformation~\eqref{three-term},
make the substitution $q\mapsto q^{2}$ and $(a,b,c,d,e)=(-q,q^2,q,-q^2,-q^3)$.
Thus
\begin{align*}
\pFq{3}{2}{-q,q^{2},q}{-q^{2},-q^{3}}{q^{2},-q}
&=\fr{(-q,-q^2;q^2)_\infty}{(-q^3,-1;q^2)_\infty} \pFq{3}{2}{q,q^{2},q}{-q^{2},-q^{2}}{q^{2},q^2} \\
&+\fr{(q,q^2,q,q^2;q^2)_\infty}{(-q^2,-q^3,-1,-q;q^2)_\infty} \pFq{3}{2}{-q,-q^{2},-q}{q^{2},-q^{2}}{q^{2},q^2} \\
&= \frac{1+q}{2}\sum_{n\geq 0} \frac{(q;q^{2})_{n}^{2}\,q^{2n}} {(-q^{2};q^{2})_{n}^{2}}
+\frac{(1+q)}{2} \fr{(q;q)_{\infty}^{2}}{(-q;q)_{\infty}^{2}}\sum_{n\ge0}
\fr{(-q;q^{2})_{n}^{2}\,q^{2n}}{(q^{2};q^{2})_{n}^{2}}.
\end{align*}
Now multiply both sides by $\fr{2q}{1+q}\fr{(-q^2;q^2)_\infty^2}{(q;q^2)_\infty^2}$ to obtain
\begin{align}
\frac{2q}{1+q}\fr{(-q^{2};q^{2})_{\infty}^{2}}
{(q;q^{2})_{\infty}^{2}}
\pFq{3}{2}{-q,\,q^{2},\,q}{-q^{2},\,-q^{3}}{q^{2};-q}
&= C_{1}(q)
+q \frac{(q^{2};q^{2})_{\infty}^{2}}{(-q;q^{2})_{\infty}^{2}}
\sum_{n\ge0}\frac{(-q;q^{2})_{n}^{2}\,q^{2n}}{(q^{2};q^{2})_{n}^{2}} \label{C1 help1}
\end{align}
Next observe that the second term on the right hand-side of~\eqref{C1 help1} simplifies as follows
\begin{align}
q \frac{(q^{2};q^{2})_{\infty}^{2}}
{(-q;q^{2})_{\infty}^{2}}
\pFq{2}{1}{-q,\,-q}{q^{2}}{q^{2};\,q^{2}}
&= q \fr{(q^{2};q^{2})_{\infty}^{2}}{(-q;q^{2})_{\infty}^{2}}
\fr{(-q,-q^3;q^{2})_{\infty}}{(q^{2};q^{2})_\infty^2}
\pFq{2}{1}{-q,q^2}{-q^3}{q^{2};-q}  \nonumber \\
&= q\sum_{n\ge0} \frac{(-q)^{n}}{1+q^{2n+1}},\label{C1 help2}
\end{align}
where in the last step we applied~\eqref{Heine-1} with $q\to q^2$ and $(a,b,c)=(-q,-q,q^2)$.
Now combining~\eqref{C1 help1} with~\eqref{C1 help2}, we deduce
\begin{equation}\label{C1 help3}
\frac{2q}{1+q}\fr{(-q^{2};q^{2})_{\infty}^{2}}
{(q;q^{2})_{\infty}^{2}}
\pFq{3}{2}{-q,q^{2},q}{-q^{2},-q^{3}}{q^{2};-q}
=C_{1}(q) +q\sum_{n\ge0} \frac{(-q)^{n}}{1+q^{2n+1}}.
\end{equation}
We now focus on the $q$-hypergeometric series
\[
\pFq{3}{2}{-q,q^{2},q}{-q^{2},-q^{3}}{q^{2};-q}
\]
which appears on the left hand-side of~\eqref{C1 help3}.
We will apply~\eqref{Watson} with $q\to q^2$ and $(a,b,c,d,e,f)=(q^2,-q^2,-q,q^{-2N},-q,q^2)$.
For simplicity, let $L$ and $R$ respectively be the left hand-side and the right hand-side of~\eqref{Watson} as with $N\to \infty$.
Then we have
\begin{align*}
L &= \lim_{N\to\infty} \pFq{8}{7}{q^2, q^3, -q^{3}, -q^2, -q, q^{-2N}, -q, q^2}{q, -q, -q^2, -q^3, q^{4+2N}, -q^3, q^2} {q^2,  -q^{2N+2}} \\
&=\lim_{N\to\infty} \sum_{n\geq 0} \fr{(q^3,-q^3;q^2)_n (-q;q^2)_n^2}{(q,-q;q^2)_n (-q^3;q^2)_n^2} (-1)^n q^{2n} (q^{-2N};q^2)_n q^{2Nn} \\
&=\lim_{N\to\infty} \sum_{n\geq 0} \fr{(1-q^{4n+2}) (1+q)^2}{(1-q^2) (1+q^{2n+1})^2}(-1)^n q^{2n} (q^{-2N};q^2)_n q^{2Nn} \\
&=\fr{(1+q)^2}{1-q^2}\sum_{n\geq 0}\fr{q^{n^2+n} (1-q^{2n+1})}{1+q^{2n+1}}.
\end{align*}
As for $R$, we have
\begin{align*}
R &= \fr{(q^4;q^2)_\infty (-q;q^2)_\infty}{(-q^3;q^2)_\infty (q^2;q^2)_\infty} \lim_{N\to\infty}\pFq{4}{3}{q,q^{-2N},-q,q^2}{-q^{2},-q^{3},-q^{1-2N}}{q^{2};q^2} \\
&= \fr{(q^4;q^2)_\infty (-q;q^2)_\infty}{(-q^3;q^2)_\infty (q^2;q^2)_\infty} \pFq{3}{2}{q,-q,q^2}{-q^{2},-q^{3}}{q^{2};-q} \\
&=\fr{1+q}{1-q^2} \pFq{3}{2}{q,-q,q^2}{-q^{2},-q^{3}}{q^{2};-q},
\end{align*}
where the second step follows since by l'Hopital's rule we have
\[
\lim_{N\to\infty}\fr{(q^{-2N};q^2)_n}{(-q^{1-2N};q^2)_n} = (-1)^n q^{-n}.
\]
Now equating $L$ and $R$ and simplifying yield
\[
(1+q)\sum_{n\geq 0}\fr{q^{n^2+n} (1-q^{2n+1})}{1+q^{2n+1}} =\pFq{3}{2}{q,-q,q^2}{-q^{2},-q^{3}}{q^{2};-q}.
\]
Finally, combining the foregoing formula with~\eqref{C1 help3} we find
\[
\Big(2q \sum_{n\geq 0}\fr{q^{n^2+n}(1-q^{2n+1})}{1+q^{2n+1}} \Big) \fr{(-q^2;q^2)_\infty^2}{(q;q^2)_\infty^2}
=C_1(q) + q\sum_{n\geq 0}\fr{(-q)^n}{1+q^{2n+1}},
\]
which is equivalent to the desired result.
\section{Proof of Theorem~\ref{thm C21}}\label{sec proof C21}
We have
\begin{align*}
C_{2}(q)
&= \sum_{n\geq 0}\fr{q^{2n+1} (-q^{2n+4};q^2)_\infty (-q^{2n+2};q^2)_\infty}{(q^{2n+1};q^2)_\infty (q^{2n+1};q^2)_\infty} \\
&=\fr{q (-q^2,-q^4;q^2)_\infty}{(q;q^2)_\infty^2} \sum_{n\geq 0} \fr{q^{2n} (q;q^2)_n^2}{(-q^2,-q^4;q^2)_n} \\
&=\fr{2q(1+q^2)(-q^2,-q^4;q^2)_\infty}{(q;q^2)_\infty^2}\sum_{n\geq 0}\fr{q^{2n} (q;q^2)_n^2}{(-1,-q^2;q^2)_{n+1}}.
\end{align*}
Then by~\eqref{Chan-Mao} applied with $q\to q^2$ and $(y,z)=(q,-1)$ and appealing to~\eqref{basic-C}, we find
\begin{align}
C_{2}(q)
&=\sum_{n\geq 0}\fr{q^{2n+1} (-q^{2n+4};q^2)_\infty (-q^{2n+2};q^2)_\infty}{(q^{2n+1};q^2)_\infty (q^{2n+1};q^2)_\infty} \nonumber \\
&=\fr{(-q^2,-q^4;q^2)_\infty}{(q;q^2)_\infty^2}\sum_{n\geq 0}\fr{q^{2n+1}(q,q;q^2)_n}{(-q^2,-q^4;q^2)_n} \nonumber \\
&=\fr{2q(-q^2;q^2)_\infty^2}{(q;q^2)_\infty^2}\sum_{n\geq 0}\fr{q^{2n}(q,q;q^2)_n}{(-1,-q^2;q^2)_{n+1}} \nonumber \\
&=\fr{2q(-q^2;q^2)^2}{(q;q^2)_\infty^2}
\Big( \fr{1}{(1+q)^2}-\fr{(q;q^2)_\infty^2}{(1+q)^2 (-1,-q^2;q^2)_\infty} \Big) \nonumber\\
&=\fr{2q(-q^2;q^2)_\infty^2}{(1+q)^2 (q;q^2)_\infty^2} - \fr{q}{(1+q)^2} \nonumber.
\end{align}
This competes the proof.
\section{Proof of Theorem~\ref{thm C31}}\label{sec proof C31}
Noting that
\begin{align*}
\sum_{n\geq 0}\fr{q^{2n+1}(-q^{2n+6},-q^{2n+2};q^2)_\infty}{(q^{2n+3},q^{2n+1};q^2)_\infty}
&= \sum_{n\geq 0}\fr{q^{2n+1}(1-q^{2n+1}) (-q^{2n+6},-q^{2n+2};q^2)_\infty}{(q^{2n+1};q^2)_\infty^2} \\
&=C_{3}(q) -\sum_{n\geq 0}\fr{q^{4n+2}(-q^{2n+6},-q^{2n+2};q^2)_\infty}{(q^{2n+1};q^2)_\infty^2},
\end{align*}
we get
\begin{equation}\label{intermed-1}
C_{3}(q)
=\sum_{n\geq 0}\fr{q^{2n+1}(-q^{2n+6},-q^{2n+2};q^2)_\infty}{(q^{2n+3},q^{2n+1};q^2)_\infty}
+\sum_{n\geq 0}\fr{q^{4n+2}(-q^{2n+6},-q^{2n+2};q^2)_\infty}{(q^{2n+1};q^2)_\infty^2}.
\end{equation}
We will now evaluate the two sums on the right hand-side of~\eqref{intermed-1}.
As for the first sum, with the help of~\eqref{Gea-Mel cmpn} applied to $q\to q^2$ and $(y,z)=(q,-q^{-2})$ and an appeal to~\eqref{basic-C}, we derive
\begin{align}
&
\sum_{n\geq 0} q^{2n+1} \fr{(-q^{2n+2},-q^{2n+6};q^2)_\infty}{(q^{2n+1},q^{2n+3};q^2)_\infty} \nonumber \\
&=\fr{q(-q^2,-q^6;q^2)_\infty}{(q,q^3;q^2)_\infty}\sum_{n\geq 0}\fr{q^{2n} (q,q^3;q^2)_n}{(-q^2,-q^6;q^2)_n} \nonumber \\
&=\fr{2(1+q^4)q(-q^2,-q^6;q^2)_\infty}{(q,q^3;q^2)_\infty} \sum_{n\geq 1}\fr{q^{2n-2} (q,q^3;q^2)_{n-1}}{(-1,-q^4;q^2)_n} \nonumber \\
&=\fr{2q^{-1}(-q^2,-q^4;q^2)_\infty}{(q,q^3;q^2)_\infty} \fr{1}{q+q^{-1}+q^2+q^{-2}}\Big(1-\fr{(q,q^3;q^2)_\infty}{(-1,-q^4;q^2)_\infty} \Big)
\nonumber \\
&=\fr{2q(-q^2,-q^4;q^2)_\infty}{(1+q)(1+q^3)(q,q^3;q^2)_\infty} -\fr{q}{(1+q)(1+q^3)} \nonumber \\
&=\fr{2q(1-q) (-q;q)_\infty^2}{(1+q)(1+q^2)(1+q^3)(-q;q^2)_\infty^2 (q;q^2)_\infty^2} -\fr{q}{(1+q)(1+q^3)} \nonumber \\
&=\fr{2q (1-q)(-q;q)_\infty (-q^4;q)_\infty}{(q^2;q^4)_\infty^2}  -\fr{q}{(1+q)(1+q^3)} \\
&=\fr{2q(-q^2;q)_\infty (-q^4;q)_\infty}{(q^2;q^4)_\infty (q^6;q^4)_\infty}-\fr{q}{(1+q)(1+q^3)} \\
&=\fr{2q}{(1+q)(1+q^3)}\fr{(-q;q)_\infty (-q^3;q)_\infty}{(q;q^2)_\infty (-q;q^2)_\infty(q^3;q^2)_\infty(-q^3;q^2)_\infty}-\fr{q}{(1+q)(1+q^3)}
\nonumber \\
&=\fr{2q}{(1+q)(1+q^3)}\fr{(-q^2;q^2)_\infty (-q^4;q^2)_\infty}{(q;q^2)_\infty(q^3;q^2)_\infty}-\fr{q}{(1+q)(1+q^3)} \label{intermed-2}.
\end{align}
As for the second sum,
using~\eqref{id appl-1} with $q\to q^2$ and $(y,z)=(q,-q^{-2})$, making an appeal to~\eqref{basic-C}, and simplifying, we obtain
\begin{align}
&\sum_{n\geq 0} q^{4n+2} \fr{(-q^{2n+2},-q^{2n+6};q^2)_\infty}{(q^{2n+1},q^{2n+1};q^2)_\infty} \nonumber \\
&=\fr{(-q^2,-q^6;q^2)_\infty}{(q,q;q^2)_\infty}\sum_{n\geq 0}\fr{q^{4n+2} (q,q;q^2)_n}{(-q^2,-q^6;q^2)_n} \nonumber \\
&=2q^2(1+q^4)\fr{(-q^2,-q^6;q^2)_\infty}{(q;q^2)_\infty^2}\sum_{n \geq 0}\fr{q^{4n}(q,q;q^2)_n}{(-1,-q^4;q^2)_{n+1}} \nonumber \\
&=2q^2(1+q^4)\fr{(-q^2,-q^6;q^2)_\infty}{(q;q^2)_\infty^2}
\Big( \fr{(1-q^2)}{(1+q)^2(1+q^3)^2}\fr{(q;q^2)_\infty^2}{(-1,-q^4;q^2)_\infty} \nonumber \\
&+\fr{(1-q^2)(1+q^{-2}) }{-q^{-2}(1+q)^2(1+q^3)^2} + \fr{1}{(1+q^3)^2} \Big) \nonumber \\
&=\fr{2q^3(2+q+q^3)}{(1+q)^2(1+q^3)^2}\fr{(-q^2,-q^4;q^2)_\infty}{(q;q^2)_\infty^2} + \fr{q^2(1-q)}{(1+q)(1+q^3)^2} \label{intermed-3}.
\end{align}
Then we achieve the desired result by adding~\eqref{intermed-2} and~\eqref{intermed-3}, and using~\eqref{intermed-1}.
\section{Remarks and conjectures}\label{sec remarks}
In this section we present several conjectures on Ramanujan-type congruences modulo $4$ (and higher powers of $2$)
for the sequences $C(k,n)$ and for the limiting sequence $c(n)$.

\begin{remark}\label{rem eta-quotient}
Write $q=e^{2\pi i\tau}$ with $\Im(\tau)>0$ and let $\eta(\tau)$ denote the Dedekind eta-function given by
\[
\eta(\tau)=q^{1/24}(q;q)_\infty
=q^{1/24}\prod_{n\ge1}(1-q^n).
\]
Using
\[
(-q^2;q^2)_\infty = \prod_{n\geq 1}(1+q^{2n}) = \frac{(q^4;q^4)_\infty}{(q^2;q^2)_\infty}
\]
and $(q;q)_\infty = (q;q^2)_\infty(q^2;q^2)_\infty$, we obtain
\[
\frac{(-q^2;q^2)_\infty^2}{(q;q^2)_\infty^2}
 = \frac{(q^4;q^4)_\infty^2}{(q;q)_\infty^2}
 = q^{-1/4}\,\frac{\eta(4\tau)^2}{\eta(\tau)^2}.
\]
In particular, up to the elementary factor $q^{-1/4}$ this common factor is an eta-quotient of weight~$0$ on $\Gamma_0(4)$. Thus suitable linear combinations of $C_1(q)$, $C_2(q)$, and $C_3(q)$ can be
expressed as modular functions on $\Gamma_0(4)$ with integral Fourier coefficients, see~Ono~\cite{Ono} for general background.
\end{remark}

{\bf 1.\ } As for $C(1,n)$, Corollary~\ref{cor 9n+8} gives the Ramanujan-type congruence
\[
C(1,9n+8) \equiv 0 \Mod{4}
\]
for all $n\geq 0$. Together with Corollary~\ref{cor C1-mod4}, this describes the possible residue classes of
$C(1,n)$ modulo $4$ in terms of the prime factorisation of $2n-1$.
%

{\bf 2.\ } Let
\[
C(q)=\lim_{k\to\infty} C_k(q) = \sum_{n\geq 0}\frac{q^{2n+1} (-q^{2n+2};q^2)_\infty}{(q^{2n+1};q^2)_\infty^2}
\]
and let $C(q) := \sum_{n\geq 0} c(n)q^n$.
Extensive computations of the coefficients $c(n)$ suggest that $c(8n+4)$ is always divisible by $4$, whereas $c(8n+6)$
appears to be divisible by $8$ (and hence also by $4$). This leads to the following Ramanujan-type conjectures for the limiting sequence $c(n)$.

\begin{conjecture}\label{conj-2}
For any nonnegative integer $n$, we have
\[
c(8n+4) \equiv 0 \Mod{4}.
\]
\end{conjecture}

\begin{conjecture}\label{conj-3}
For any nonnegative integer $n$, we have
\[
c(8n+6) \equiv 0 \Mod{8}.
\]
\end{conjecture}
These congruences are similar to Ramanujan's classical congruences for the ordinary partition function, such as
\[
p(5n+4)\equiv 0 \pmod{5},\qquad p(7n+5)\equiv 0 \pmod{7},
\]
but here they concern the limiting two-color partition sequence  $c(n)$. It would be interesting to obtain a direct modular-forms proof of Conjectures~\ref{conj-2} and~\ref{conj-3}.
\medskip

It is worth mentioning that $c(8n+4)$ is not always divisible by $8$. For instance, for $n=1$ we have $8n+4=12$ and
\[
c(12)=284\equiv 4 \pmod{8},
\]
so the two conjectures above appear to be genuinely distinct.

{\bf 3.\ } Note that we may write
\[
C(q) = \fr{q(-q^2;q^2)_\infty}{(q;q^2)_\infty^2} \sum_{n\geq 0} \fr{q^{2n}(q;q^2)_n^2}{(-q^2;q^2)_n} = A(q) S(q),
\]
where
\begin{align*}
A(q)&= \fr{q(-q^2;q^2)_\infty}{(q;q^2)_\infty^2} = \sum_{n\geq 0} a(n) q^n, \\
S(q) &= \sum_{n\geq 0}\fr{q^{2n}(q;q^2)_n^2}{(-q^2;q^2)_n} = \sum_{n\geq 0} s(n) q^n.
\end{align*}
Then one can show that for any $n\geq 0$,
\[
a(4n) \equiv a(4n+3) \equiv 0\Mod{4}
\]
As for $S(q)$, we have the following conjecture.
\begin{conjecture}\label{conj-4}
For any nonnegative integer $n$, we have
\[
s(4n+1) =0.
\]
\end{conjecture}
One can show that
\begin{equation}\label{Sq-id}
S(q)= \fr{(q;q^2)_\infty^2}{(q^2;q^2)_\infty^2}\sum_{n\geq 0}\fr{1+q^{2n+1}}{1-q^{2n+1}} q^{3n^2+3n}\sum_{j=-n}^n (-1)^j q^{-2j^2}.
\end{equation}
\begin{remark}
The double sum on the right-hand side of~\eqref{Sq-id} can be seen as a Hecke-type series where the exponent
on q is an indefinite quadratic form
\[
Q(n,j)=3n^2+3n-2j^2.
\]
It is therefore natural to expect modular or mock-modular behaviour of $S(q)$, and in particular that the congruences in Conjectures~\ref{conj-2}--\ref{conj-5} might be dealt with using the arithmetic of the quadratic form $Q$.
\end{remark}

We close with the following conjecture about the even parts of $(q^4;q^4)_\infty S(q)$.
\begin{conjecture}\label{conj-5}
For any nonnegative integer $N$, we have
\begin{align*}
[q^{2N}] \Big((q^4;q^4)_\infty S(q) \Big) &= [q^{2N}] \Big( \sum_{n\geq 0}(-1)^n q^{6n^2+4n}\sum_{j=-n}^n (-1)^j q^{-2j^2} \\
& \qquad\qquad - \sum_{n\geq 1}(-1)^n q^{6n^2-4n}\sum_{j=-n+1}^{n-1} (-1)^j q^{-2j^2} \Big) \\
&=[q^{2N}] \Big( \sum_{n\geq 0}(-1)^n q^{6n^2+4n}(1+q^{4n+2})\sum_{j=-n}^n (-1)^j q^{-2j^2} \Big).
\end{align*}
\end{conjecture}
\bigskip


\noindent{\bf Data Availability Statement.\ }
Not applicable.

\bigskip

\noindent{\bf Declarations.\ }
The authors state that there is no conflict of interest.
\end{document}